# Selmer Groups of Hecke Characters and Chow Groups of Self Products of CM Elliptic Curves [*]

Jonathan Dee

January 4, 1999

## Contents

## 1 Introduction

Formulae relating special values of $L$-fuctions to orders of Selmer groups have been in existence for about a century. Dirichlet's formula can be recast in this form, as can the Birch and Swinnerton-Dyer conjecture from the 1960's. Conjectural generalisations of these formulae were made by Deligne, Beilinson, then Bloch and Kato, extending the framework to the case of a general mixed motive over a number field. They are beautiful in their own right, but recently novel applications have begun to become apparent. The formula predicted by Bloch and Kato in the case of an adjoint representation is very important in deformation theory, and for the symmetric square of an elliptic curve appears in Wiles's work [Wiles]. In this paper we shall give application to the study of Chow groups, after proving various formulae for Selmer groups attached to elliptic crurves with complex multiplication.

---

[*]*Mathematics Subject Classification* 11G40(14C25)

Let $E$ be an elliptic curve with complex multiplication by $K$, a quadratic imaginary extension of $\mathbb{Q}$, and suppose that $E$ is defined over $K$. If $L(E/K, 1)$ is nonzero then it was shown in [CW] that $E(K)$ is finite, and later Rubin showed in [Rubin] that the whole Selmer group is finite with order as expected by the Birch and Swinnerton-Dyer conjecture (up to a controlled constant).

On the other hand, if $I_{\mathfrak{f}}$ denotes the group of fractional ideals of $K$ prime to $\mathfrak{f}$, the conductor of $E$, then there is a certain homomorphism

$$\psi : I_{\mathfrak{f}} \longrightarrow K^{\times},$$

the Grössencharacter attached to $E$. Hecke attached an $L$-function $L(\psi, s)$ to $\psi$ and proved both its analytic continuation and functional equation. It is a result of Deuring that there is an equality of $L$-functions

$$L(\psi, s)L(\overline{\psi}, s) = L(E/K, s). \tag{1}$$

One may also consider $L$-functions attached to the homomorphisms $\psi^k$ for positive integers $k$ and ask whether analogues of the Birch and Swinnerton - Dyer conjecture hold. This takes us into the domain of very general conjectures of Beilinson, and Bloch and Kato (see [BK] and [FPR]).

Our main result on $L$-functions is the following.

**Theorem 1.0.1**
*Let $k$ be a positive integer and $p > k$ a prime where $E$ has good supersingular reduction. Then*

$$\#H^1_f(K, A_k) = \left( \frac{L(\overline{\psi}^k, k)L(\psi^k, k)}{\overline{\Omega}^k \Omega^k} \right)_p,$$

*with the understanding that both sides may be simultaneously infinite.*

Here $H^1_f(K, A_k)$ is the Selmer group defined in Section 3.3 and $\Omega^k$ is certain period of $E$. See Section 2 for the notation $(x)_p$ for $x \in \overline{\mathbb{Q}}_p$. It is not obvious that the expression in brackets is a rational integer, but this fact follows from Damerell's Theorem (see Corollary 7.18 in [Rubin2]).

We further show that both sides of the above equality are finite if $k \geq 3$, and extend this to the case $k = 2$ when $E$ is assumed to be defined over $\mathbb{Q}$. In the final two sections we bound the Selmer group of the dual representation by using the functional equation, then deduce results over $\mathbb{Q}$ when $E$ is assumed to be defined over $\mathbb{Q}$.

The method we use is modelled on Rubin's approach to the Birch and Swinnerton - Dyer conjecture. When Rubin deals with the case of supersingular primes he makes use of Wiles's explicit reciprocity law. This does not apply to the Galois representations we are interested in, so we must use Kato's generalisation of Wiles's result. This and the two variable Main Conjecture proved by Rubin are the key ingredients in our proof.

It is essential for us to use some $p$-adic Hodge theory - the correct Selmer group at a supersingular prime cannot even be defined without some use of $B_{cris}$ and $B_{dR}$, and Kato's reciprocity law makes use of the theory in its statement. The analogous results to ours when $E$ has good ordinary reduction at $p$ was proved by Li Guo in [LiGuo]. In this case the local condition at $p$ in the definition of the Selmer group may be given without resorting to $p$-adic Hodge theory, and the proof is consequently less technically involved.

After the first part of this work (up to and including Section 4) was completed it came to the attention of the author that some similar results have been obtained independently by Han (see [Han]). We subsequently realised that his results together with results of Nekovář [N2] could be used to provide information about Chow groups of arbitrary codimension for self products of CM elliptic curves. Our main result in this direction is the following. See Chapter 5 for more introductory material.

**Theorem 1.0.2**
*Let $E$ be as above, with good reduction at a prime $p > 3$, and let $X$ be the $d$-fold self product of $E$. Let $i$ be a strictly positive integer. Assume that for all $0 \leq n \leq i - 1$*

- *$p^2 - 1 \nmid n + 1 + p(n - 1)$ if $E$ has supersingular reduction at $p$,*

- *$p - 1 \nmid 2n$ if $E$ has ordinary reduction at $p$.*

*(e.g. if $p > 2i + 1$).*
*Then*
$$\mathrm{Im}(cl_X(CH^i(X)) \cap F^2 H^{2i}(X_{et}, \mathbb{Q}_p(i)) = 0,$$
*where $X$ is the $d$-fold product of $E$. The same is true for $E \otimes_{\mathbb{Q}} K$.*

*Furthermore, if $p > 2d + 1$ then in $H^{2i}(\overline{X}, \mathbb{Z}_p(i))$ we have that $\mathrm{Im}(cl_X)_{tors}$ is finite, and zero for all*
$$p \nmid (2\pi)^j \frac{L(\psi^{2j}, -j)}{\Omega^{2j}}.$$

Here the filtration is induced by the Hoschild-Serre spectral sequence, and

$$cl_X : CH^i(X) \longrightarrow H^{2i}(\overline{X}_{et}, \mathbb{Z}_p(i)) \longrightarrow H^{2i}(\overline{X}_{et}, \mathbb{Q}_p(i))$$

is the étale cycle class map.

I would like to thank Jan Nekovář for his encouragement, frequent help, and tireless patience, and John Coates for very helpful comments.

## 2  Notations and Conventions

We fix embeddings of $\overline{\mathbb{Q}}$ into both $\overline{\mathbb{Q}}_p$ and $\mathbb{C}$.

If $p$ is a rational prime and $x \in \overline{\mathbb{Q}}_p$, then write

$$(x)_p = p^{v_p(x)},$$

where $v_p$ is the $p$-adic valuation of $\overline{\mathbb{Q}}_p$ normalised so that $v_p(p) = 1$, with the convention that $(0)_p = \infty$.

If $K$ is a number field and $v$ is a prime of $K$, denote by $K_v$ the completion of $K$ at $v$. If $F$ is a local or global field let $\mathcal{O}_F$ denote the ring of integers of $F$. If $F$ is any field we shall write $\overline{F}$ for some algebraic closure of $F$. If $K/F$ is a Galois extension, write $G(K/F)$ for the Galois group. Write $G_F = G(\overline{F}/F)$ whenever $F$ is a perfect field.

We use the notation $\mathbb{F}_q$ for the finite field with $q$ elements.

As usual any Galois cohomology groups will always be taken to mean continuous cohomology groups. We normalise our Hodge-Tate weights so that for example $\mathbb{C}_p(1)$ has Hodge-Tate weight 1 (rather than $-1$).

We apologise for using both the arithmetic and geometric Frobenius. This is because the geometric Frobenius is related to the $f$ operator on $D_{cris}$) (this being related, of course, to the Frobenius on the reduction of the variety mod $p$ via the crystalline cohomology), while the arithmetic Frobenius sits more naturally with the Grössencharacter.

# 3  Background

## 3.1  $p$-adic Theory

We refer the reader to [Fon1] for the construction and properties of Fontaine's rings of $p$-adic periods $B_{cris}$ and $B_{dR}$, and to [FPR] for the basics of their application to local Galois cohomology.

Let $K$ be a finite extension of $\mathbb{Q}_p$, and let $K_0$ be the maximal absolutely unramified subextension of $\mathbb{Q}_p$ in $K$ with Frobenius $\sigma$ (lifting $x \mapsto x^p$).

If $R$ is any topological ring then by an $R$-representation of $G_K$ we mean an $R$-module $V$ of finite type equipped with a continuous $R$-linear action of $G_K$. We shall sometimes refer to a $\mathbb{Q}_p$-representation as a $p$-adic representation.

We first generalise Proposition 3.8 of [BK] to the case of representations with coefficients. Let $F$ be a finite extension of $\mathbb{Q}_p$. Suppose that $V$ is an $F$-representation of $G_K$, and let $T \subset V$ be an $\mathcal{O}_F$-lattice in $V$, such that the action of $G_K$ on $V$ makes $T$ into a $\mathcal{O}_F$-representation of $G_K$. Write $A = V/T$, a discrete $\mathcal{O}_F$-module with linear action of $G_K$. We have the following:

**Proposition 3.1.1**
*Suppose that $F$ is an unramified extension of $K$. Then the Tate pairing with coefficients*

$$H^1(K, A) \times H^1(K, T^*(1)) \longrightarrow F/\mathcal{O}_F, \tag{2}$$

*where $T^* = Hom_{\mathcal{O}_F}(T, \mathcal{O}_F)$, is perfect and induces an isomorphism*

$$\frac{H^1(K, A)}{H^1_e(K, A)} \overset{\sim}{\longrightarrow} Hom_{\mathcal{O}_F}(H^1_g(K, T^*(1)), F/\mathcal{O}_F). \tag{3}$$

*Proof:* First note that

$$Hom_{\mathcal{O}_K}(\mathcal{O}_F, \mathcal{O}_K) \cong \mathcal{O}_F. \tag{4}$$

Indeed, given $x \in \mathcal{O}_F$, consider the function given by

$$y \mapsto Tr_{F/K}(xyc),$$

where $c$ generates the inverse different of $F/K$. By definition of the different and trace this map is well defined and $\mathcal{O}_K$-linear. To prove that all maps as above are of this form, note that trace gives a perfect $K$-linear pairing on $F$, so any function as above is of the form

$$y \mapsto Tr_{F/K}(zy)$$

for some $z \in F$. By the definition of the inverse different $z = cx$ for some $x \in \mathcal{O}_F$.

Next let $M$ be any $\mathcal{O}_F$-module of finite rank. The adjunction formula tells us

$$Hom_{\mathcal{O}_F}(M, Hom_{\mathcal{O}_K}(\mathcal{O}_F, \mathcal{O}_K)) \cong Hom_{\mathcal{O}_K}(M, \mathcal{O}_K), \qquad (5)$$

where, for instance, if

$$g : M \longrightarrow \mathcal{O}_K$$

is $\mathcal{O}_K$-linear, then we may define, for each $m \in M$

$$f(t) : s \mapsto g(st),$$

$t \mapsto f(t)$ being the desired $\mathcal{O}_F$-linear map.

Substituting $T$ for $M$ into (4) and (5), and noting that both equalities will be Galois equivariant if we let Galois act on $F$ and $K$ trivially, we find that

$$Hom_{\mathcal{O}_F}(T, \mathcal{O}_F) \cong Hom_{\mathcal{O}_K}(T, \mathcal{O}_K).$$

Similarly

$$Hom_{\mathcal{O}_F}(T, F/\mathcal{O}_F) \cong Hom_{\mathcal{O}_K}(T, K/\mathcal{O}_K).$$

The result now follows formally from the classical Tate duality pairing, and the result of Bloch and Kato. $\qquad \square$

See [BK] for the definitions of the exponential and dual exponential maps for $p$-adic representations

$$exp_V : \frac{D_{dR}(V)}{D_{dR}^0(V)} \longrightarrow H_e^1(K, V),$$

$$exp_V^* : \frac{H^1(K, V)}{H_g^1(K, V)} \longrightarrow D_{dR}^0(V).$$

They are natural in the variable $V$, so respect the action of $F$ if $V$ is an $F$-representation for $F$ as above. One has the following formulae (see for example [N1]), whereall dimensions are taken over $\mathbb{Q}_p$.

**Lemma 3.1.2**

$\dim H_e^1(K, V) = \dim(D_{dR}(V)/D_{dR}^0(V)) + \dim H^0(K, V) - \dim D_{cris}(V)^{f=1}$,

$\dim H_f^1(K, V) = \dim(D_{dR}(V)/D_{dR}^0(V)) + \dim H^0(K, V)$,

$\dim H_g^1(K, V) = \dim(D_{dR}(V)/D_{dR}^0(V)) + \dim H^0(K, V) + \dim D_{cris}(V^*(1))^{f=1}$.

**Remark 3.1.3** Counting dimensions in the exact sequence of $\mathbb{Q}_p$-vector spaces

$$0 \longrightarrow D_{cris}(V)^{f=1} \longrightarrow D_{cris}(V) \longrightarrow (1-f)D_{cris}(V) \longrightarrow 0$$

shows that $\dim D_{cris}(V)^{f=1} = \dim D_{cris}(V)/(1-f)D_{cris}(V)$.

We shall need to understand the characteristic polynomial of the crystalline Frobenius acting on $D_{cris}(V)$. Let $X$ be a smooth projective variety defined over $K$ and supppose $X$ has a smooth projective model over $\mathcal{O}_K$. If $l$ is any prime write $V_{l,n} = H^n_{et}(X_{\overline{K}}, \mathbb{Q}_l)$. Write $H^n_{cris}(X)$ for the $n^{th}$ rational crystalline cohomology group of the special fibre of a smooth projective model of $X$ over $\mathcal{O}_K$. It is a $K_0$-vector space with a $\sigma$-semilinear Frobenius operator $f$. Thus if $a = [K_0 : \mathbb{Q}_p]$, then $f^a$ acts $\sigma$-linearly. Since $X$ has good reduction the $G_K$-representation $V_l$ is unramified, so if $Frob$ denotes a *geometric* Frobenius in $G_{\mathbb{Q}_p}$, then there is a well defined action of $Frob^a$ on $V_{l,n}$. Katz-Messing [KM] deduced the following from the Weil conjectures.

**Fact 3.1.4** If $l \neq p$

$$\mathrm{char}(1 - Tf^a \mid H^n_{cris}(X)) = \mathrm{char}(1 - TFrob^a \mid V_{l,n}).$$

Combining this with the Crystalline Conjecture proved by Faltings (and Fontaine - Messing, Kato, Tsuji) that $D_{cris}(V_p)$ is naturally isomorphic to $H_{cris}(X)$, in a maner compatible with the action of Frobenius, we deduce

**Lemma 3.1.5** *We have*

$$\mathrm{char}(1 - Tf^a \mid D_{cris}(V_p)) = \mathrm{char}(1 - TFrob^a \mid V_{l,n}).$$

$\square$

Deligne's proof of the Weil conjectures tells us that the roots of these characteristic polynomials are algebraic integers, all of whose conugates have absolute value equal to $p^{-an/2}$.

In fact we shall only use Lemma 3.1.5 when $V_{l,n}$ is the rational Tate module of an elliptic curve with complex multiplication, and in this special case it can be deduced in an ad hoc way from the results of Fontaine in [Fon3].

## 3.2 Grössencharacters and $L$-functions

For background information on elliptic curves with complex multiplication see chapter 2 of [Sil2].

From now on let $K$ be a quadratic imaginary extension of $\mathbb{Q}$ of discriminant $d_K$, and let $E$ be an elliptic curve over $K$ of conductor $\mathfrak{f}$ with complex multiplication by the maximal order in $K$. The field $K$ then has class number 1. Fix an embedding of $K$ into $\mathbb{C}$, and a minimal Weierstrass equation for $E$. This gives rise to a lattice $L$ in $\mathbb{C}$ and uniformisation

$$\mathbb{C}/L \cong E(\mathbb{C}). \tag{6}$$

By the theory of complex multiplication we may choose $\Omega \in \mathbb{C}$ such that $L = \Omega \mathcal{O}_K$. The induced action of $End_K(E)$ on the tangent space of $E$ induces an identification of $K$ with $End_K(E) \otimes_{\mathbb{Z}} \mathbb{Q}$ which we suppress from the notation.

Attached to $E$ is a Grössencharacter

$$\psi : I_{\mathfrak{f}} \longrightarrow K^{\times} \subset \mathbb{C},$$

such that for every prime ideal $\mathfrak{p}$ of $K$, $\psi(\mathfrak{p})$ is a generator of $\mathfrak{p}$ and $\psi$ gives the action of the *arithmetic* Frobenius at $\mathfrak{p}$ on the Tate module of $E$ at primes $l$ with $\mathfrak{p} \nmid l$.

Hecke proved analytic continuation and the functional equation for $L$-functions attached to Grössencharacters. If $\Lambda$ denotes the extended $L$-function then

$$\Lambda(\overline{\psi}^k, s) = \epsilon(\overline{\psi}, s)\Lambda(\psi^k, k+1-s), \tag{7}$$

where

$$\epsilon(\overline{\psi}^k, s) = W(\overline{\psi}^k)(\mid d_K \mid \cdot N_{K/\mathbb{Q}}\mathfrak{f})^{1/2-s}.$$

Here $d_K$ is the discriminant of $K$, and $W(\overline{\psi}^k)$ is the root number which in our situation is $\pm 1$.

The connection between the $L$-functions of $\psi^k$ and $\overline{\psi}^k$ is

$$\overline{L(\psi^k, \overline{s})} = L(\overline{\psi}^k, s). \tag{8}$$

If $E$ happens to be defined over $\mathbb{Q}$ then by considering conjugation of Frobenius elements by complex conjugation one shows that $\psi(\overline{\mathfrak{p}}) = \overline{\psi}(\mathfrak{p})$, Hence for $s$ giving

$$L(\psi, s) = L(\overline{\psi}, s). \tag{9}$$

Combining (8) with (9) we deduce that if $E$ is defined over $\mathbb{Q}$ then $L(\psi^k, k)$ is real.

## 3.3 The Galois representations

From now on we suppose that $E$ has good supersingular reduction at the prime $\mathfrak{p}$ of $K$ with $\mathfrak{p} \nmid 6$. It is known that this is equivalent to assuming that $\mathfrak{p}$ is inert: $\mathfrak{p} = p\mathcal{O}_K$ for some rational prime $p$. For simplicity we identify $p$ with $\mathfrak{p}$. It will be clear in any given context whether $p$ means a rational prime or the prime in $\mathcal{O}_K$ which it generates. Recall that another equivalent condition for an elliptic curve to have supersingular reduction at a prime $\mathfrak{p}$ above $p$ is that the Tate module of $p$-power torsion points on the reduced curve should vanish. Write $\mathcal{O}_P$ for the ring of integers in $K_p$.

We shall be interested in various representations related to the Tate module of $E$, $T_E$. Let $T_E$ denote the $\mathbb{Z}_p$-adic Tate module of $E$, let $V_E = T_E \otimes_{\mathbb{Z}_p} \mathbb{Q}_p$, and define $V_k = V_E^{\otimes k}$, a free module over $K_p$ of rank 1 with continuous action of $G_K$, the tensor product being taken over $K_p$. Similarly define $T_k = T_E^{\otimes k}$, the tensor product being over $\mathcal{O}_p$, and $A_k = V_k/T_k$.

We have fixed an embedding of $\overline{\mathbb{Q}}$ into $\overline{\mathbb{Q}}_p$, and this gives an injection $G_{\mathbb{Q}_p} \hookrightarrow G_{\mathbb{Q}}$. In particular we may consider the restriction of the $G_K$-representations $T_k$, $V_k$ and $A_k$ to $G_{K_p}$. Since $E$ has good reduction at $p$ the $G_{K_p}$-representation $V_E$ is crystalline. It follows from [Fon2] that $V_k$ and $V_k^*(1)$ are also crystalline.

Let $V$ be an arbitrary $p$-adic representation of $G_F$, where $F$ is a number field. Suppose we are given a Galois equivariant $\mathbb{Z}_p$-lattice $T$ in $V$, and put $A = V/T$. Define the Selmer group for $A$ over $K$ to be

$$H_f^1(F, A) = \ker\left(H^1(F, A) \longrightarrow \prod_{v \text{ prime of } F} \frac{H^1(F_v, A)}{L_v}\right),$$

where $L_v = i_v(W_v)$,

$$i_v : H^1(F_v, V) \longrightarrow H^1(F_v, A),$$

and

$$W_v = \begin{cases} \ker(H^1(F_v, V) \longrightarrow H^1(I_v, V)) & v \nmid p \\ H_f^1(F_v, V) & v = p \end{cases},$$

$I_v$ being the inertia subgroup of $G_{F_v}$.

We shall also need the restricted Selmer group. This is given by

$$H_f^1(F, A)' = \ker\left(H^1(F, A) \longrightarrow \prod_{v \nmid p} \frac{H^1(F_v, A)}{L_v}\right).$$

The connection between these two Selmer groups may then be summarised by the exact sequence

$$0 \longrightarrow H_f^1(F, A) \longrightarrow H_f^1(F, A)' \longrightarrow \bigoplus_{v \mid p} \frac{H^1(F_p, A)}{H_f^1(F_p, A)}. \tag{10}$$

In Section 5 we shall be interested in the Galois representations $W_k = V_{2k}(1-k)$ and $W_k^*(1)$ for positive integers $k$. We now intend to calculate the dimensions of the local cohomology groups

$$H_x^1(K_p, V_k), \qquad H_x^1(K_p, V_k^*(1))$$

$$H_x^1(K_p, W_k), \qquad H_x^1(K_p, W_k^*(1))$$

for $x = e, f$ or $g$. Since we shall have to prove results which are sometimes identical for $V_k$ and $W_k$ and sometimes slightly different we shall write $U$ for either $W_k$ or $V_k$, and similarly $T_U$ and $A_U$ for their canonical lattices and co-lattices. The first step is to understand $D_{cris}(U)$ and $D_{cris}(U^*(1))$ over $K_p$. Both $K_p$-modules are in fact $K_p \otimes_{\mathbb{Q}_p} K_p$-modules, where the first copy of $K_p$ is acting by functoriality (from the complex multiplication of $E$), and the second copy is $K_p = B_{cris}^{G_{K_p}}$, which acts on $(B_{cris} \otimes_{\mathbb{Q}_p} U)^{G_{K_p}}$ in the natural way. Since the Frobenius $f$ is linear with respect to the first $K_p$, anti-linear with respect to the second, and $\dim_{\mathbb{Q}_p} D_{cris}(U) = 4$, we deduce that $D_{cris}(U)$ and $D_{cris}(U^*(1))$ are free $K_p \otimes_{\mathbb{Q}_p} K_p$-modules of rank 1. Since $(\psi(p)) = (p)$ we have $\psi(p) = up$ for some unit $u$ in $\mathcal{O}_K$. From section 3.2 we know that the *geometric* Frobenius class $Frob_p$ in $G_{K_p}$ acts on the rational $l$-adic Tate-module of $E$ $V_{E,l}$ by $\psi^{-1}$ (for $l \neq p$), so by Lemma 3.1.5 we find

$$\text{char}(1 - Tf^2 \mid D_{cris}(V_k)) = (1 - (up)^{-k}T)^2, \tag{11}$$

Now, if $E$ is defined over $\mathbb{Q}$ then it is known that $u = -1$. We shall only need this fact when $k = 2$, so we make the following

**Hypothesis 3.3.1** If $k = 2$ assume that $E$ is defined over $\mathbb{Q}$.

Thus $f^2$ acts on $D_{cris}(V)$ by $(up)^{-k}$. Using the Weil pairing we have

$$V_k^*(1) \cong V_k(1 - k), \tag{12}$$

so $f^2$ acts on $D_{cris}(V^*(1))$ by $u^k p^{2-k}$. It furthermore follows that $f$ acts on $D_{cris}(V_1)$ with respect to a suitable $K_p \otimes K_p$-basis by

$$\begin{pmatrix} 0 & 1 \\ -1/p & 0 \end{pmatrix},$$

and hence $f$ acts on $D_{cris}(W_k^*(1))$ by $(-1)^k\sigma$, a fact which will be needed in Section 5.

By (11) and the above determination of the structure of $D_{cris}(U)$ as $K_p \otimes K_p$-module equipped with the action of Frobenius we deduce from the above discussion:

**Proposition 3.3.2**

$$D_{cris}(V_k)^{f=1} = D_{cris}(W_k)^{f=1} = 0,$$
$$D_{cris}(V_k^*(1))^{f=1} = 0 \text{ for } k = 1 \text{ and } k \geq 3,$$
$$\dim_{K_p} D_{cris}(W_k^*(1))^{f=1} = 1.$$

We next examine $H^0(K_p, U)$ and $H^0(K_p, U^*(1))$. We shall prove that both groups are trivial for all $k$.

**Lemma 3.3.3** $H^0(K_p, U) = H^0(K_p, U^*(1)) = 0$ *for all* $k \geq 1$.

*Proof:* The hypothesis of complex multiplication implies that the groups are either trivial or equal to the whole space. We are thus reduced to proving that $G_{K_p}$ acts nontrivially on $U$ and $U^*(1)$. Were the action trivial then the Hodge-Tate weights of the representation would both equal zero, which is not the case. $\square$

Kato proves in section 2.2.3 of [Kato] that

$$D_{dR}^0(V_k) \cong (Cotan(\hat{E}) \otimes_{\mathbb{Z}_p} \mathbb{Q}_p)^{\otimes_{K_p} k}$$

canonically, where $Cotan(\hat{E})$ is the cotangent space for the formal group of $E/K$ at the prime $p$, a free $\mathcal{O}_p$-module of rank 1. hence

$$\dim_{K_p} D_{dR}^0(V) = 1. \tag{13}$$

Since $V_k$ and $V_k^*(1)$ are de Rham representations of $G_{K_p}$

$$\dim_{K_p} D_{dR}(V_k) = \dim_{K_p} D_{dR}(V_k^*(1)) = 2. \tag{14}$$

We could also see the above from the fact that $V$ is de Rham with Hodge-Tate weights $0$ and $k$, but we shall need the above isomorphism anyway in the construction of the explicit reciprocity map. The representation $W_k$ is de Rham with Hodge-Tate weights $-k$ and $k$, so both $D_{dR}^0(W_k)$ and $D_{dR}^0(W_k^*(1))$ also have $K_p$-dimension equal to $1$.

We are now in a position to prove:

**Proposition 3.3.4**
*In the following all dimensions are with respect to the $K_p$-vector space structures.*

$$\dim H^1(K_p, U) = \dim H^1(K_p, U^*(1)) = 2.$$

*If $k \neq 2$ and $x = e, f,$ or $g$ then*

$$\dim H_x^1(K_p, V_k) = \dim H_x^1(K_p, V_k^*(1)) = 2.$$

$$\dim H_x^1(K_p, W_k) = \begin{cases} 1 & x = e, \text{ or } f, \\ 2 & x = g. \end{cases}$$

$$\dim H_x^1(K_p, W_k^*(1)) = \begin{cases} 0 & x = e, \\ 1 & x = f, \text{ or } g. \end{cases}$$

*Proof:* Tate's calculation of the local Euler-Poincaré characteristic and local duality gives

$$\dim H^1(K_p, U) = 2 + \dim H^0(K_p, U) + \dim H^0(K_p, U^*(1)) = \dim H^1(K_p, U^*(1)),$$

and both $H^0$'s vanish by Lemma 3.3.3. the rest of the formulae follow immediately by substituting (13), Proposition 3.3.2, and Lemma 3.3.3 into the formulae of Proposition 3.1.2. $\qquad\square$

### 3.4 Iwasawa Theory

Let $K_\infty = K(A)$ be the extension of $K$ given by the kernel of

$$\psi^k : G_K \longrightarrow \mathcal{O}_p^\times, \tag{15}$$

(which is of course also the kernel of $\overline{\psi}^k$). Hence

$$G_\infty = G(K_\infty/K) \subset \mathcal{O}_p^\times$$

canonically. We know from the theory of complex multiplication that so long as $p$ is prime to $6\mathfrak{f}$, $\psi$ and $\overline{\psi}$ are actually isomorphisms (see Corollary 5.17 of [Rubin2]). We may thus decompose

$$G_\infty = \Delta \times \Gamma,$$

where $\Delta$ is a finite abelian group of order prime to $p$, and $\Gamma$ is (noncanonically) isomorphic to $\mathbb{Z}_p^2$ (recall $k > 0$).

Let $K_0$ denote the fixed field of $\Gamma \subset G_\infty$, so $G(K_\infty/K_0) \cong \Gamma$, and more generally write $K_n$ for the fixed field of $\Gamma^{p^n}$. It is known that there is a unique (totally ramified) prime $p_n$ of $\mathcal{O}_{K_n}$ above $p$. Write $K_{n,p}$ for the completion of $K_n$ at that prime.

Write $\Lambda = \mathbb{Z}_p[[G_\infty]]$ for the Iwasawa algebra corresponding to $G_\infty$. Since $\#\Delta$ is prime to $p$, $\mathbb{Z}_p[\Delta]$ is semi-simple. Let $e_\chi$ be the idempotent in $\mathbb{Z}_p[\delta]$ corresponding to the character $\chi$ of an irreducible $\mathbb{Z}_p$-representation of $\Delta$. Then

$$\Lambda = \mathbb{Z}_p[\Delta][[\Gamma]] = \bigoplus_\chi \Lambda_\chi,$$

where $\chi$ runs over all characters of irreducible $\mathbb{Z}_p$ representations of $\Delta$ and

$$\Lambda_\chi = e_\chi \Lambda = R_\chi[[\Gamma]],$$

$R_\chi$ being either $\mathbb{Z}_p$ or the ring of integers in the unique unramified quadratic imaginary extension of $\mathbb{Z}_p$ according as $\dim(\chi) = 1$ or $2$. More generally, we have a decomposition for any $\Delta$-module $Y$:

$$Y = \bigoplus Y^\chi, \tag{16}$$

where $\chi$ runs over all irreducible $\mathbb{Z}_p$-representations of $\Delta$ and $Y = e_\chi Y$.

Write $\mathcal{A}_k$ for the kernel of the homomorphism

$$\psi^k : \Lambda_{\chi_k} \longrightarrow \mathcal{O}_p$$

deduced from (15). Here $\chi_k$ is the restriction of $\psi^k$ to $\Delta$.

We shall restrict attention from now on to values of $k$ such $\chi_k$ is $\mathbb{Z}_p$-irreducible. This is not too drastic because of the following.

**Lemma 3.4.1** *$\chi_k$ is irreducible as a $\mathbb{Z}_p$-representation of $\Delta$ if and only if $p+1 \nmid k$.*

*Proof:* This follows easily from the fact that $\psi$ is surjective, $E[p] \cong \mathcal{O}/p$, and that $Aut(E[p]) \cong (\mathcal{O}/p)^\times$ is a cyclic group of order $p^2 - 1$.    □

The following two simple results will be needed later.

**Lemma 3.4.2** *Identify $E[p]$ with $\mathcal{O}/p$. Then $\chi_k = \overline{\chi}\chi$ if and only if $p^2 - 1 \mid k - 1 - p$.*

*Proof:* Suppose $\chi_k = \overline{\psi}\psi$. We know that $\psi$ maps $G_K$ onto $\Delta = (\mathcal{O}/p)^*$. Now complex conjugation induces the unique nontrivial continuous automorphism of $K_p$ which is the Frobenius for $K_p/\mathbb{Q}_p$. Hence if $x \in \Delta$, $c(x) = x^p$. By the surjectivity of $\chi_1$, $\chi_k = \chi_1\overline{\chi_1}$ if and only if $x^{k-1} = x^p$ for all $x$, if and only if $k - 1 \equiv p \bmod p^2 - 1$.

The converse is clear by the same reasoning.    □

**Lemma 3.4.3** *Suppose $p > k$ and $p \nmid 6$. Then the homomorphism*

$$\psi^k : \Lambda \longrightarrow \mathcal{O}_p$$

*is surjective.*

*Proof:* This follows from the fact that $\psi : G_K \to \mathcal{O}_p^\times$ is surjective, as remarked at the beginning of this section. Hence if $p \nmid k$, then $1 + p\mathcal{O}_p$ is contained in the image of our map, and if the torsion subgroup of $\mathcal{O}_p$ is not contained in the image, then $\chi_k$ maps $G_K$ into $\mathbb{F}_p^\times$. We saw in the proof of Proposition 3.4.1 that this cannot occur if $p > k - 1$. □

**Hypothesis 3.4.4** Here we collect our assumptions on the prime $p$. In view of Lemmas 3.4.1 and 3.4.2 we suppose from now on that $p > k$. Recall that we are already assuming that $p \nmid 6$ is a prime of good supersingular reduction for $E$.

### Some $\Lambda$-modules

We shall use use the following standard $\Lambda$-modules for the $\Delta \times \mathcal{O}_K$-extension $K_\infty/K$. See [Rubin] for further explanation.

- $A_\infty$ ($p$-parts of ideal class groups).

  It is known that $A_\infty$ is a finitely generated torsion $\Lambda$-module.

- $U_\infty$ (principal local units).

  It shown in [Win] that $U_\infty$ is finitely generated and torsion free as a $\Lambda$-module, and that $U_\infty^\chi$ is free of rank 2 over $\Lambda_\chi$. (Note that in order to apply Wintenberger's results we use Lemma 3.4.2 which ensures that $\chi$ is not the cyclotomic character).

- $E_\infty$ (closure of global units).

  Note that since the local units at level $n$, $U_n$ split as the product of principal units and roots of unity of order prime to $p$, the image $E_n$ of the global units in $U_n$ may be considered equally as a subgroup or quotient of the global units.

- $X_\infty$.

  This $\Lambda$-module is finitely generated of rank 1, and has no nonzero pseudo-null submodules (see [Green]).

We shall also need the Iwasawa module of elliptic units. For more details on the construction and properties of these units see [Rubin2].

Firstly, if $x$ is the coordinate function on $E$ coming from the Weierstrass model fixed in Section 3.2, define for each ideal $\mathfrak{a}$ of $K$ prime to 6 ($\mathfrak{a}$ will always be such an ideal in this section)

$$\Theta_{E,\mathfrak{a}} = \alpha^{-12} \Delta(E)^{N\mathfrak{a}-1} \prod_{P \in E[\mathfrak{a}]-O} (x - x(P))^{-6},$$

where $(\alpha) = \mathfrak{a}$.

$\Theta_{E,\mathfrak{a}}$ is independent of $\alpha$ since all units in $\mathcal{O}_K$ are killed by raising to the twelfth power. In fact in general it is independent of the choice of the Weierstrass model and defined over any field of definition of $E$.

Fix a set $B$ of ideals of $K$ prime to $\mathfrak{af}$ such that the Artin map induces a bijection of $B$ with $G(K(\mathfrak{f})/K)$. Let $P_{\mathfrak{f}}$ generate $E[\mathfrak{f}]$ as $\mathcal{O}_K$-module, and define

$$\Lambda_{E,\mathfrak{a}}(P) = \prod_{\mathfrak{b} \in B} \Theta_{E,\mathfrak{a}}(\psi(\mathfrak{b})P_{\mathfrak{f}} + P).$$

This formula gives a rational function, and one can check that it is in fact defined over $K$. It is independent of the choice of $B$, (but not necessarily of $P_{\mathfrak{f}}$), and if $P_{\mathfrak{b}} \notin E[\mathfrak{f}]$ (where $\mathfrak{b}$ is as above), then $\Lambda_{E,\mathfrak{a}}(P_{\mathfrak{b}})$ is a global unit in $K(E[\mathfrak{b}])$.

Recall that

$$\xi : \mathbb{C}/L \cong E(\mathbb{C})$$

is our fixed uniformisation of $E(\mathbb{C})$. It is compatible with the group structure. Since $\psi(p^n)$ generates the ideal $(p^n)$,

$$P_n = \xi(\psi(p^n)^{-1}\Omega) \tag{17}$$

will have exact order $p^n$. One may show that $\eta_{\mathfrak{a},n} = \Lambda_{E,\mathfrak{a}}(P_n)$ is a global unit in $K_n$, and moreover that for fixed $\mathfrak{a}$ the $\eta_{\mathfrak{a},n}$ are norm compatible in the tower $K_\infty/K$. Let $C'_n$ denote the group of units in $K_n$ generated as $\mathbb{Z}[G(K_n/K)]$-module by the roots of unity in $K_n$ and the $\eta_{\mathfrak{a},n}$, where $\mathfrak{a}$ runs over the ideals of $K$ prime to $6p\mathfrak{f}$. Finally define $C_n$ to be the closure of the image of $C'_n$ in $U_n$, and

$$C_\infty = \varprojlim C_n.$$

We have the following result giving the structure of $C_\infty$ as a $\Lambda_\chi$-module. First we need a preliminary result.

**Lemma 3.4.5** *There exists an ideal $\mathfrak{a}$ prime to $6\mathfrak{f}p$ such that $N\mathfrak{a} - \psi(\mathfrak{a})^k$ is prime to $p$.*

*Proof:* It is clear that $N\mathfrak{a}$ is congruent to $\psi^k(\mathfrak{a})$ modulo $p$ for all $\mathfrak{a}$ if and only if $\chi_k$ is the cyclotomic character. By the Čebotarev density theorem it is sufficient to check this for primes coprime to $6\mathfrak{f}p$ since there are only a finite number of exceptions. Hence the lemma follows from Lemma 3.4.2 (and Hypothesis 3.4.4). □

**Lemma 3.4.6** *$C_\infty^\chi$ is free of rank 1. It is generated by $(\eta_{\mathfrak{a},n})$ where $\mathfrak{a}$ is any ideal of $\mathcal{O}_K$ prime to $6p\mathfrak{f}$ such that*

$$\psi^k(\mathfrak{a}) \not\equiv N_{K/\mathbb{Q}}(\mathfrak{a}) \bmod p.$$

*Proof:* One may show that for all ideals $\mathfrak{a}$ and $\mathfrak{b}$ of $\mathcal{O}_K$ not dividing $6p\mathfrak{f}$ and for all n the following relation holds

$$(\eta_{\mathfrak{a},n})^{Frob_{\mathfrak{b}} - N\mathfrak{b}} = (\eta_{\mathfrak{b},n})^{Frob_{\mathfrak{a}} - N\mathfrak{a}} \tag{18}$$

(see [Rubin2] Theorem 7.4 and Lemma 7.10). On the other hand, by Lemma 3.4.5 we may choose $\mathfrak{a}$ to satisfy $\psi^k(\mathfrak{a}) \not\equiv N\mathfrak{a} \bmod p$. This implies that the image of $Frob_{\mathfrak{a}} - N\mathfrak{a}$ is invertible in $\Lambda_\chi$, and hence (by (18) that $C_\infty^\chi$ is generated as $\Lambda_\chi$-module by $\eta_{\mathfrak{a},n}$ and the image of the $p$-power roots of unity in $K_n$. Since $\chi$ is not the cyclotomic character by Lemma 3.4.2 we deduce that $C_n^\chi$ is generated by $\eta_{\mathfrak{a},n}$. The norm compatibility of the $\eta_{\mathfrak{a},n}$ shows that the $\eta_{\mathfrak{a},n}$ fit together to give an element $\eta_{\mathfrak{a}}$ of $C_\infty$. A simple compactness argument then shows that the image of $\eta_{\mathfrak{a}}$ generates $C_\infty^\chi$.

Finally we have $C_\infty^\chi \subset U_\infty^\chi$, the latter being a free $\Lambda_\chi$-module. Hence $C_\infty^\chi$ is torsion free and the proof is complete. $\qquad\square$

## 4   Bounding Selmer Groups

### 4.1   Kato's Reciprocity Map

Write $P = (P_n)_n$ for the generator of the $\mathcal{O}_p$-module $T_E$ defined in (17), and from now on fix $\chi = \chi_k$. In [Kato] Kato considers a homomorphism

$$l : U_\infty \longrightarrow Cotan(\hat{E})^{\otimes k} \otimes K_p,$$

(the tensor products being taken over $\mathcal{O}_p$), as the composite of

$$h : U_\infty \overset{\alpha}{\longrightarrow} \varprojlim H^1(K_{p,n}, \mathbb{Z}_p(1)/p^n)$$
$$\overset{\beta}{\longrightarrow} \varprojlim H^1(K_{p,n}, T_k^*(1)/p^n) \overset{\gamma}{\longrightarrow} \varprojlim H^1(K_p, T_k^*(1)/p^n) \overset{\sim}{\longrightarrow} H^1(K_p, T_k^*(1))$$

and

$$exp^* : H^1(K_p, T_k^*(1)) \longrightarrow H^1(K_p, V_k^*(1)) \longrightarrow D_{dR}^0(V_k^*(1)) \overset{\sim}{\longrightarrow} Cotan(\hat{E})^{\otimes k} \otimes K_p.$$

where

- the inverse limits are with respect to corestriction and the natural projections.

- $\alpha$ is deduced from the Kummer maps for the fields $K_{p,n}$, $\beta$ is the inverse limit of cup product with $P_n^{\otimes(-k)}$, and $\gamma$ is corestriction.

One knows that $Cotan(\hat{E})$ is an invertible $\mathcal{O}_p$-module, hence $Colie(G)^{\otimes k} \otimes K_p$ is a free $K_p$-module of rank 1.

We shall show that the map

$$h_\chi : U_\infty^\chi \longrightarrow H^1(K_p, T_k^*(1))$$

induced from $h$ is surjective (we shall also exhibit its kernel) and thus obtain the equality

$$h_\chi(\ker(l_\chi)) = \ker(exp^*) = H_g^1(K_p, T_k^*(1)), \tag{19}$$

where $l_\chi$ is the restriction of $l$ to $U_\infty^\chi$.

First note that restriction induces an isomorphism

$$H^1(K_p, A_k) \overset{\sim}{\longrightarrow} H^1(K_{\infty,p}, A_k)^{G_\infty} \overset{\sim}{\longrightarrow} Hom_{\mathbb{Z}_p}(G_{K_{\infty,p}}^{ab}, A_k)^{G_\infty}.$$

This follows from the proof of Lemma 6.2 of [Rubin2] and the following lemma .

**Lemma 4.1.1** *Let $\Delta \subset (\mathcal{O}/p^n)^\times$ be the image of*

$$G_K \xrightarrow{\psi^k} \mathcal{O}_p^\times \longrightarrow (\mathcal{O}_p/p^n)^\times.$$

*Let $\Delta$ act on $\mathcal{O}_p/p^n$ by*

$$z : x \mapsto z^k x.$$

*Then (under Hypothesis 3.4.4) $H^1(\Delta, \mathcal{O}_p/p^n) = 0$.*

*Proof:* Fist note that $\Delta$ is not a $p$-group by the proof of Lemma 3.4.1. Let $\Delta'$ be the prime to $p$-part of $\Delta$. Then $\Delta'$ is isomorphic to the image of $\chi$. Since $\mathcal{O}_p/p^n$ is a $p$-group, $H^i(\Delta', \mathcal{O}_p/p^n) = 0$ for $i \geq 1$. Furthermore, $H^0(\Delta', \mathcal{O}_p/p^n) = 0$. Indeed, suppose otherwise. Then $\Delta'$ must be killed by $k$. The order of $\Delta'$ is a multiple of $p + 1$, so $(p+1) \mid k$. This contradicts Hypothesis 3.4.4.

The inflation/restriction sequence

$$0 \longrightarrow H^1(\Delta/\Delta', (\mathcal{O}_p/p^n)^{\Delta'}) \longrightarrow H^1(\Delta, \mathcal{O}_p/p^n) \longrightarrow H^1(\Delta', \mathcal{O}_p/p^n)$$

now gives us the desired result. $\qquad\qquad\qquad\qquad\qquad\qquad\qquad\qquad\qquad\square$

By local class field theory

$$G_{K_{\infty,p}}^{ab} \cong U_\infty,$$

and this isomorphism is $G_\infty$-equivariant. Hence

$$H^1(K_p, A_k) \cong Hom_{\mathbb{Z}_p}(U_\infty, A_k)^{G_\infty} \cong Hom_{\mathcal{O}_p}(U_\infty^\chi/\mathcal{A}_k U_\infty^\chi, A_k),$$

where, we recall,

$$\mathcal{A}_k = \ker(\psi^k : \Lambda_\chi \longrightarrow \mathcal{O}_p).$$

We deduce an isomorphism of $\mathcal{O}_p$-modules

$$U_\infty^\chi/\mathcal{A}_k \cong Hom_{\mathcal{O}_p}(H^1(K_p, A_k), K_p/\mathcal{O}_p),$$

and hence after applying local Tate duality an isomorphism

$$t : U_\infty^\chi/\mathcal{A}_k U_\infty^\chi \xrightarrow{\sim} H^1(K_p, T_k^*(1)).$$

**Proposition 4.1.2**
*The following diagram commutes:*

$$
\begin{array}{ccc}
U_\infty^\chi & \xrightarrow{\;\;h_\chi\;\;} & H^1(K_p, T_k^*(1)) \\
\downarrow & & \text{id} \downarrow \\
U_\infty^\chi/\mathcal{A}_k & \xrightarrow[\;\sim\;]{\;\;t\;\;} & H^1(K_p, T_k^*(1)).
\end{array}
$$

*Proof:* Given the definition of Tate duality it suffices to prove

$$c \cup h(u) = c^{res}(u) P^{\otimes -k} \qquad \text{for all } c \in H^1(K_p, A_k) \text{ and } u \in U_\infty^\chi \qquad (20)$$

where $c^{res}$ is the image of $c$ in $Hom_{\mathbb{Z}_p}(U_\infty, A_k)^{G_\infty}$. Note that $P$ occurs in the definition of $h$, so it is not surprising that it appears in the right hand side of (20).

Let $(g(u)_n)_n = \beta\alpha(u)_n$, where $\alpha$ and $\beta$ are the maps occuring in the definition of $h$. Writing $cores_n$ (resp. $res_n$) for corestriction (resp. restriction) on cohomology from $K_{n,p}$ to $K_p$, the projection formula for Galois cohomology tells us that for any $n$:

$$c \cup h(u) = c \cup cores_n(g(u)_n) = cores_n(res_n(c) \cup g(u)_n). \qquad (21)$$

Let $c'_n = res_n(c) \cup P_n^{\otimes -k}$, so $res_n(c) \cup g(u)_n = c'_n \cup \alpha(u)_n$.

By the continuity of $c$ we may find $n_0 \geq 1$ such that $c \in \mathrm{im}(H^1(K_p, T_k/p^{n_0}))$. From now on let $n \geq n_0$. Thus $c \in H^1(K_p, T/p^n)$.

Lemma 4.1.3 below will tell us that $c'_n \cup \alpha(u)_n = c'_n(u_n)$. Substituting this into (21) we find

$$c \cup h(u) = c'_n(u_n) = c^{res}(u)P^{\otimes -k},$$

and the proof of the proposition is complete.

$\square$

**Lemma 4.1.3** *Let $F$ be a finite extension of $\mathbb{Q}_p$, and suppose $[F(\mu_{p^n}) : F]$ is prime to $p$. Then the pairing*

$$\Phi : H^1(F, \mathbb{Z}/p^n) \times F^* \longrightarrow \mathbb{Z}/p^n\mathbb{Z}$$

$$(c, z) \longmapsto inv(c \cup \delta(z))$$

*is just $\Phi(c, z) = c(\rho(z))$, where $\rho : F^* \longrightarrow G_F^{ab}$ is the local reciprocity map, $\delta$ is the Kummer homomorphism, and we identify $H^1(F, \mathbb{Z}/p^n)$ with $Hom_{\mathbb{Z}_p}(G_F^{ab}, \mathbb{Z}/p^n)$.*

*Proof:* $F(\mu_{p^n})/F$ is Galois, and $H = G(F(\mu_{p^n})/F)$ has order prime to $p$ by hypothesis, so

$$H^1(F, \mathbb{Z}/p^n) = H^1(F(\mu_{p^n}), \mathbb{Z}/p^n)^H.$$

Now, a property of the invariant map and restriction implies

$$\Phi(c, z) = (\#H)^{-1}inv(res(c) \cup res(\delta(z))),$$

(this makes sense since $\#H$ is prime to $p$), and the result follows from Proposition 5 in Chapter XIV of [Serre].

$\square$

Proposition 4.1.2 immediately implies (19). Applying Tate local duality, and using Proposition 3.3.4, we deduce

**Proposition 4.1.4**

$$H_f^1(K_p, A_k) = \ker(H^1(K_p, A_k) \longrightarrow Hom_{\mathbb{Z}_p}(\ker(l), A_k)).$$

$\square$

### 4.2   The Selmer group and Kato's explicit reciprocity law

First note that the choice of Weirstrass model we have made gives us an identification of $\mathcal{O}(\hat{E})$ with $\mathcal{O}_p[[w]]$, where $w = x/y$ is a local parameter at the identity of $E$. Since $E$ has supersingular reduction at $p$ we may identify $T_E$ with the Tate module of $\hat{E}$ as representation of $G_{K_p}$.

Let $\Lambda_{p,\mathfrak{a}} \in K_p((w))$ denote the expansion of $\Lambda_{\mathfrak{a}}$ in terms of $w$. One may show explicitly that $\Lambda_{p,\mathfrak{a}} \in \mathcal{O}_p[[w]]^\times$. With $\eta_{n,\mathfrak{a}}$ defined as in Section 3.4, define $\eta_{\mathfrak{a}} = (\eta_{n,\mathfrak{a}})_n$ It follows from the definition of $\Lambda_{\mathfrak{a}}$ that $\Lambda_{p,\mathfrak{a}}$ is the Coleman power series of $\eta_{\mathfrak{a}}$ relative to $(P_n)$.

Next we recall the defintion of the Coates-Wiles homomorphisms. Let $D = \frac{1}{\lambda'(w)} \frac{d}{dw}$ be the logarithmic derivative operator on $\mathcal{O}_{K_p}[[w]]$, where

$$\lambda(w) \in w + w^2 K_p[[w]]$$

is logarithm map of $\hat{E}$. In fact

$$\lambda'(w) \in 1 + w \mathcal{O}_p[[w]],$$

so it is clear that $D$ does indeed map $\mathcal{O}_p[[w]]$ into itself. The $k^{th}$ Coates-Wiles homomorphism is defined for $k \geq 1$ by

$$\delta : U_\infty \longrightarrow \mathcal{O}_p \tag{22}$$

$$u \mapsto \left( D^k \log(g_u(w)) \right) \mid_{w=0}.$$

(We suppress $k$ from the notation since it is fixed). Here $\log g_u(w)$ is meant in a formal sense: since $k \geq 1$,

$$D^k \log g_u(w) = D^{k-1} \left( \frac{1}{\lambda'(w)} \frac{d}{dw} g_u \right)(w).$$

The Coates-Wiles homomorphisms are compatible with the action of $\Lambda$ in the following way:

$$\delta(fu) = \psi^k(f)\delta(u) \tag{23}$$

for all $f \in \Lambda$. Using Lemma 3.4.3 and (23) we find:

**Lemma 4.2.1** *If $u$ is any element in $U_\infty$ then $\delta(\Lambda_\chi u)$ is an ideal in $\mathcal{O}_p$.*                    □

The following explicit reciprocity law of Kato gives a cohomological interpretation to $\delta$. See [Kato] for the proof.

**Theorem 4.2.2**
*If*

$$l : U_\infty \longrightarrow Cotan(\hat{E})^{\otimes k} \otimes K_p$$

*is the homomorphism defined in section 4.1, then*

$$l(u) = \frac{1}{(k-1)!} \omega^{\otimes k} \otimes \left( \left( \frac{d}{\omega} \right)^k \log(g_u) \right)(O_E),$$

*where $\omega$ is any nonzero element of $Cotan(\hat{E})$, $O_E$ is the identity of the group law on $E$, and $g_u$ denotes the Coleman power series of $u$ relative to $P_n$*                    □

Now, if $\omega_E$ is the canonical translation invariant holomorphic differential attached to $E/K$ by our choice of Weierstrass equation, then writing $\omega_E$ in terms of $dx$ and $dy$, and expanding $x$ and $y$ in power series in $w$, one finds that

$$\omega_E = \lambda'(w)dw$$

(see page 8 of [Rubin2]), and this in fact characterises $\lambda$ as an element of $z + z^2 K_p[[w]]$. Since $\lambda'(w) \in 1 + z\mathcal{O}_p[[w]]$, and by its derivation from $\omega_E$, we have that $\lambda'(w)dw \in Cotan(\hat{E})$.

We may thus identify $Cotan(\hat{E})^{\otimes k}$ with $\mathcal{O}_p$ by choosing $(\lambda'(w)dw)^{\otimes k}$ as a basis. Under this identification the operator $\frac{d}{\omega}$ in Theorem 4.2.2 becomes $D$. Evaluating elements of $\mathcal{O}(\hat{E})$ at the identity correponds to evaluating power series in $w$ at $w = 0$ under $\mathcal{O}(\hat{E}) \cong \mathcal{O}_p[[w]]$ since $w$ is a local parameter at the identity. We may thus rephrase the conclusion of Theorem 4.2.2 as

$$l(u) = \frac{1}{(k-1)!}\delta(u). \tag{24}$$

**Remark 4.2.3** Notice that we deduce from this that the image of $l$ is a subset of $\mathcal{O}_p$, a fact which is by no means obvious from the cohomological definition of $l$.

We next need to identify the image of $l$. This is achieved by the following proposition which was shown to us by Coates.

**Proposition 4.2.4**
If $k < p^2 - 1$ then the homomorphism

$$l : U_\infty \longrightarrow K_p$$

has image equal to $\mathcal{O}_p$.

*Proof:* We show that $\delta$ has image $(k-1)!\mathcal{O}_p$ then appeal to (24).

The tower $K_{p,n}$ is totally ramified and abelian, so by Lubin-Tate theory $K_{p,n}$ corresponds to adjoining to $K_p$ the $p^n$ torsion of a Lubin-Tate formal group $\mathcal{F}$ whose Frobenius lifting is given by the polynomial $\pi X + X^{p^2}$, for some uniformizer $\pi$ of $K_p$. By the uniqueness property of the Coleman power series, the construction of the Coates-Wiles homomorphism and Theorem 4.2.2 we see that $\delta$ coincides with the Coates-Wiles homomorphism derived from $\mathcal{F}$ (see [deShalit]) with respect to the generator $v = (v_n)$ corresponding to $P$ under the given isomorphism $\hat{E} \cong \mathcal{F}$.

Let $\beta$ be the unique $(p^2 - 1)^{th}$ root of $1 - \pi$ satisfying $\beta \equiv 1 \bmod p$, and let $v = (v_n)$ denote a generator of the Tate module of $\mathcal{F}$. The $v_n$ is a uniformiser of $K_{p,n}$ by Lubin-Tate theory. A simple calculation shows that

$$N_{n,n-1}(\beta - v_n) = \beta - v_{n+1}.$$

Thus, using the fact that $v_n$ is a uniformizer of $K_{p,n}$, we see $u = (\beta - v_n) \in U_\infty$. Visibly $g_u(w) = \beta - w$ is the Coleman power series for $u$ relative to $\mathcal{F}$.

If $w = \phi(z)$ gives the exponential map for $\mathcal{F}$, so $\phi = \lambda^{-1}$, then

$$\frac{d}{dz}\log g_u(w) = -\frac{1}{\lambda'(w)}\frac{-1}{\beta - w} = -\frac{1}{\beta}\sum_{r=0}^{\infty}\frac{w^r}{\beta^r},$$

since $1/\lambda'(w) = \phi'(z)$. On the other hand, we have $[\zeta]_{\mathcal{F}}(w) = \zeta w$ for all $\zeta \in \mu_{p^2-1}$, so $\phi(\zeta z) = \zeta\phi(z)$, for all $z$. Comparing coefficients one both sides of this equality yields

$$\phi(z) = z + \sum_{n=1}^{\infty} a_n z^{1+n(p^2-1)}$$

for some $a_n$. Thus $\frac{d}{dz}\log g_u(w)$ has coefficients of $w^r$ and $z^r$ equal for $1 \le r < p^2-1$, and $\delta(u) = (k-1)! \cdot \text{unit}$. The result now follows immediately, since by Lemma 4.2.1 the image of $\delta$ must be some ideal in $\mathcal{O}_p$. $\qquad\square$

We next wish to relate $H^1_f(K, A_k)$ to an Iwasawa module. Firstly, the same proof as in Lemma 6.4 of [Rubin2] tells us that

$$H^1_f(K, A_k)' \cong Hom_{\mathbb{Z}_p}(X_\infty, A_k)^{G_\infty},$$

in view of Lemma 4.1.1.

We now show the existence of a suitable quotient of $X_\infty$, with $\Lambda$-rank 1.

**Lemma 4.2.5** *There exists a decomposition $U^\chi_\infty = V_1 \oplus V_2$ where $V_1$ and $V_2$ are free $\Lambda^\chi$-modules of rank 1, $l(V_2) = 0$, and $E^\chi_\infty \not\subset V_2$. Furthermore $\ker(l_\chi) = \mathcal{A}_k V_1 + V_2$.*

*Proof:* In section 3.4 we remarked that $U^\chi_\infty$ is free of rank 2 over $\Lambda_\chi$. Fix a splitting

$$U^\chi_\infty = \Lambda_\chi v_1 \oplus \Lambda_\chi v_2.$$

It follows from Lemma 4.2.1 and (24) that $l(\Lambda_\chi v_1)$ and $l(\Lambda_\chi v_2)$ are ideals in $\mathcal{O}_p$. Proposition 4.2.4 tells us that $l$ is surjective, so we must have either $l(\Lambda_\chi v_1)$ or $l(\Lambda_\chi v_2)$ equal to the whole of $\mathcal{O}_p$. Suppose the former, without loss of generality.

We may choose $f \in \Lambda_\chi$ such that $l(v_2) = l(f v_1)$. Then if

$$f_0 \in \ker(\psi^k : \Lambda_\chi \longrightarrow \mathcal{O}_p)$$

it follows from (23) and (24) that we have $l(V_2) = l((f + f_0)v_1)$. We may choose $f_0$ so that, writing $g = f + f_0$, we have $E^\chi_\infty \not\subset \Lambda_\chi(v_2 - gv_1)$. The proof is now finished, if we take $V_1 = \Lambda_\chi v_1$ and $V_2 = \Lambda_\chi(v_2 - gv_1)$. $\qquad\square$

Fix a decomposition of $U^\chi_\infty$ as in lemma 4.2.5, and let $\tilde{U} = U^\chi_\infty/V_2$. Class field theory gives us a homomorphism of $\Lambda$-modules

$$U^\chi_\infty \longrightarrow X^\chi_\infty,$$

so we may define $\tilde{X} = X^\chi_\infty/\text{im}(V_2)$. Since $E^\chi_\infty \not\subset V_2$ we deduce that $\text{Im}(V_2)$ has rank 1, and that $\tilde{X}$ is a torsion $\Lambda_\chi$-module. Also, by definition

$$X^\chi_\infty/(\mathcal{A}_k + \text{im}(\ker(l_\chi))) = \tilde{X}/\mathcal{A}_k.$$

By (10):

$$H^1_f(K, A_k) = \ker\left(Hom_{\mathbb{Z}_p}(X_\infty, A_k)^{G_\infty} \longrightarrow \frac{H^1(K_p, A_k)}{H^1_f(K_p, A_k)}\right).$$

Proposition 4.1.4 yields

$$
\begin{aligned}
H_f^1(K, A_k) &= \ker(Hom_{\mathbb{Z}_p}(X_\infty, A_k)^{G_\infty} \longrightarrow Hom_{\mathbb{Z}_p}(\ker(l), A_k)) \\
&= Hom_{\mathbb{Z}_p}(X_\infty^\chi/(\mathcal{A}_{\|k} + \operatorname{im}(\ker(l))), A) \\
&= Hom_{\mathbb{Z}_p}(\tilde{X}, A_k)^\Gamma.
\end{aligned}
\tag{25}
$$

We next use the two-variable main conjecture as proved by Rubin to relate the Selmer group to the quotient $U_\infty^\chi/E_\infty^\chi$. See [Rubin] and [Rubin2].

**Theorem 4.2.6 (The Main Conjecture)**

$$
\operatorname{char}(A_\infty) = \operatorname{char}(E_\infty/C_\infty).
$$

$\square$

A standard exact sequence in global class field theory runs

$$
0 \longrightarrow E_\infty/C_\infty \longrightarrow U_\infty/C_\infty \longrightarrow X_\infty \longrightarrow A_\infty \longrightarrow 0,
$$

(see [Rubin2] page 40). We immediately deduce from this an exact sequence

$$
0 \longrightarrow E_\infty^\chi/C_\infty^\chi \longrightarrow \tilde{U}/\operatorname{im}(C_\infty^\chi) \longrightarrow \tilde{X} \longrightarrow A_\infty^\chi \longrightarrow 0.
\tag{26}
$$

Indeed the only question is exactness at $E_\infty^\chi/C_\infty^\chi$. This follows from the choice of $V_1$ and $V_2$: since $E_\infty^\chi$ is torsion-free of rank 1 (it is a submodule of $U_\infty^\chi$ which is free), and $E_\infty^\chi \not\subset V_2$ we must have either that $E_\infty^\chi/E_\infty^\chi \cap V_2$ is torsion or $E_\infty^\chi \cap V_2 = 0$. But $U_\infty^\chi/V_2$ is torsion-free, so the latter possibility must hold, and exactness follows. From Theorem 4.2.6 and (26) we deduce

$$
\operatorname{char}_{\Lambda_\chi}(\tilde{X}) = \operatorname{char}_{\Lambda_\chi}(\tilde{U}/\operatorname{im}(C_\infty^\chi)).
\tag{27}
$$

We shall use this equality to deduce that

$$
\#H_f^1(K, A_k) = \#Hom(\tilde{U}/\operatorname{im}(C_\infty^\chi), A_k)^\Gamma.
\tag{28}
$$

More generally, let $K_0 \subset L \subset K_\infty$ be a subextension such that $\Gamma_L = G(L/K_0) \cong \mathbb{Z}_p^d$ for $d = 1$ or 2. Write $\Lambda_L = \mathbb{Z}_p[[G(L/K)]]$.

Given two such subextensions of $K_\infty$ with $F \subset L$ and $G(L/F) \cong \mathbb{Z}_p$, the natural projection $G(L/K) \longrightarrow G(F/K)$ induces a surjective homomorphism of $\mathbb{Z}_p$-algebras

$$
\pi_{L/F} : \Lambda_L \longrightarrow \Lambda_F.
$$

Consider the following condition for torsion $\Lambda_L$-modules $X$:

$(pseud(L, X))$     $X$ has no pseudo-null $\Lambda_L$-submodules.

The following is due to Perrin-Riou (see [PR] and [Rubin] Lemma 6.2).

**Proposition 4.2.7**
*Suppose that $M$ is a finitely-generated torsion $\Lambda$-module such that $pseud(L, M)$ holds whenever $L$ is a $\mathbb{Z}_p^d$-extension of $K_0$ in $K_\infty$, ($d = 1$ or 2). Then*

$$
\#Hom(M, A_k)^{G_\infty} = \#(\mathcal{O}_\mathfrak{p}/(\overline{\psi}^k(\operatorname{char}_{\Lambda_\chi}(M^\chi)))),
$$

*where both sides are possibly infinite.* $\square$

**Corollary 4.2.8**

$$\#Hom(\tilde{X}, A_k)^\Gamma = \#Hom(\tilde{U}_\infty/im(\mathcal{C}_\infty), A_k)^\Gamma.$$

*Proof:* By (27) and Proposition 4.2.7 it is enough to show that $pseud(L, M)$ holds for $M = \tilde{X}$ or $\tilde{U}_\infty/\mathrm{im}(C_\infty^\chi)$, and $L$ a $\mathbb{Z}_p^d$ extension of $K_0$ in $K_\infty$, $d = 1$ or 2. Fix such an $L$. Now $pseud(L, X_\infty)$ holds for $d = 2$ by fundamental results of Greenberg (see [Green]), and $V_2$ is free, so $pseud(K_\infty, \tilde{X})$ holds. Furthermore $pseud(L, \tilde{X})$ holds for $d = 1$ since $(X_\infty)_L \cong X_\infty(L)$ as $\Lambda_L$ modules, where $X_\infty(L)$ denotes the maximal $p$-unramified abelian $p$-extension of $L$. Finally, $pseud(L, \tilde{U}/\mathrm{im}(C_\infty^\chi))$ holds since $\tilde{U}_\infty$ and $C_\infty^\chi$ are both free $\Lambda^\chi$-modules. $\qquad\square$

We may now deduce (28), as promised, from Corollary 4.2.8 and (25). On the other hand $l$ induces a homomorphism

$$l : \tilde{U} \longrightarrow \mathcal{O}_p.$$

Since $l(U_\infty) = \mathcal{O}_p$ and $l(V_2) = 0$ it follows that $l(\tilde{U}) = \mathcal{O}_p$. Now $l$ induces

$$l : \tilde{U}/\mathcal{A}_k \hookrightarrow \mathcal{O}_p,$$

so this injection is an isomorphism. Hence (28) implies that

$$\#H^1_f(K, A_k) = [\mathcal{O}_p : l(\mathcal{C}_\infty^\chi)] \qquad (29)$$

## 4.3 Connection with $L$-values

It remains to compare the image of the elliptic units under $l$ with $L(\overline{\psi}^k, k)$. This is achieved by a calculation due to Coates and Wiles. The following is Theorem 7.22 of [Rubin2], in light of the reciprocity law and the fact that $\Lambda_{a,\mathfrak{a}}$ is the Coleman power series of $\eta_\mathfrak{a}$.

**Theorem 4.3.1**
We have

$$l(\eta_\mathfrak{a}) = 12(-1)^{k-1} f^k (N\mathfrak{a} - \psi(\mathfrak{a})^k) \Omega^{-k} L(\overline{\psi}^k, k),$$

where $f$ is such that $P_\mathfrak{f} = \xi(\Omega/f)$. $\qquad\square$

**Corollary 4.3.2**

$$l(C_\infty^\chi) = (\Omega^{-k} L(\overline{\psi}^k, k)) \mathcal{O}_p.$$

*Proof:* First note that we may ignore the $p$-power roots of unity which occur in the definition of $C_n$. Indeed, $l_n$ is a homomorphism, factors through $C_n^\chi$, and since $\chi$ is not the cyclotomic character by Proposition 3.4.2 it is trivial on roots of unity. A compactness argument then shows that $\mathcal{C}_\infty$ is generated by the images of $\eta_\mathfrak{a}$ where $\mathfrak{a}$ runs over ideals prime to $6\mathfrak{f}p$. The image of $\mathcal{C}_\infty$ under $l$ is thus the ideal generated by the image of the $\eta_\mathfrak{a}$. We are thus finished by Theorem 4.3.1 (noting that $p$ is prime to $6\mathfrak{f}$) and Lemma 3.4.5.

$\qquad\square$

We thus conclude our main result on special values of $L$-functions.

**Theorem 4.3.3**
For $p$ satisfying Hypothesis 3.4.4

$$\#H^1_f(K, A_k) = \left(\frac{L(\overline{\psi}^k, k)L(\psi^k, k)}{\overline{\Omega}^k\Omega^k}\right)_p,$$

with the understanding that both sides may be simultaneously infinite.

*Proof:* Substituting the result of Corollary 4.3.2 into (29) we have

$$\#H^1_f(K, A_k) = \left[\mathcal{O}_p : \left(\frac{L(\overline{\psi}^k)}{\Omega^k}\right)\mathcal{O}_p\right].\tag{30}$$

Since $L(\overline{\psi}^k, k) = \overline{L(\psi^k, k)}$ by (8) we have

$$\#H^1_f(K, A_k) = \left[\mathcal{O}_p : \frac{L(\overline{\psi}^k, k)}{\Omega^k}\mathcal{O}_p\right]$$
$$= \left(\frac{L(\overline{\psi}^k, k)L(\psi^k, k)}{\overline{\Omega}^k\Omega^k}\right)_p.$$

$\square$

**Corollary 4.3.4** *If $k \geq 2$ (and Hypothesis 3.4.4 holds) then $H^1_f(K, A_k)$ is always finite, and its order is equal to $\left(\frac{L(\overline{\psi}^k, k)L(\psi^k, k)}{\overline{\Omega}^k\Omega^k}\right)_p$.*

*Proof:* The Euler product for $L(\overline{\psi}^k, s)$ converges absolutely for $s > k/2 + 1$ yielding the assertion for $k > 2$. Furthermore it is a classical result following from Tate's thesis (see [Lang] page 313) that $L(\phi, n)$ is finite and nonzero whenever $\phi$ is a nontrivial Gr"ossencharacter of weight $2n$ and $n$ is an integer. This yields the assertion when $k = 2$. $\square$

**Remark 4.3.5** For $k = 1$ the Selmer group $H^1_f(K, A_1)$ will be infinite if the Mordell-Weil group of $E/K$ has positive rank, and the converse is true if the Tate-Shafarevich group of $E/K$ has finite $p$-part.

## 4.4  Duality

At the beginning of this section we assume only that $p$ is a prime where $E$ has good supersingular reduction. As in earlier sections we deal with the case of $V_k$ and the weight $-2$ representations denoted $W_k$ simultaneously. Let $U$ denote any of these representations.

**Proposition 4.4.1**
For all $k \geq 1$ we have

$$\dim H^1_f(K, U) = \dim H^1_f(K, U^*(1))$$

*Proof:* We use the Poitou-Tate global duality sequence

$$0 \longrightarrow H^1_f(K, U) \longrightarrow H^1(G_S, U) \longrightarrow \bigoplus_{v \in S} \frac{H^1(K_v, U)}{H^1_f(K_v, U)}$$

$$\longrightarrow H^1_f(K, U^*(1)) \longrightarrow H^2(G_S, U) \longrightarrow \bigoplus_{v \in S} H^2(K_v, U) \longrightarrow 0. \quad (31)$$

Here $S$ is the finite set of primes consisting of primes where $U$ is ramified and $p$. Thus if $v \in S$ and $v \neq p$, then $U^{I_v} \neq U$. But $U$ is free of rank 1 over $K_p$ and $U^{I_v}$ is a $K_p$-subspace of $U$. This implies that $U^{I_v} = 0$, and

$$H^1_f(K_v, U) = H^1(G_{K_v}/I_v, U^{I_v}) = 0.$$

On the other hand, by Tate's calculation of the local Euler-Poincaré characteristic in the unequal residue characteristic case we have

$$\dim H^1(K_v, U) = \dim H^2(K_v, U), \quad (32)$$

while from Proposition 3.1.2 we have

$$\dim \frac{H^1(K_p, U)}{H^1_f(K_p, U)} = \dim H^1_f(K_p, U^*(1)) = 2 + \dim H^2(K_p, U). \quad (33)$$

Recall that our convention is to take all dimensions over $\mathbb{Q}_p$, unless we say otherwise. Equations (32) and (33) together imply

$$\dim \left( \bigoplus_{v \in S} \frac{H^1(K_v, U)}{H^1_f(K_v, U)} \right) = 2 + \left( \bigoplus_{v \in S} H^2(K_v, U) \right). \quad (34)$$

Next, Tate's calculation of the global Euler-Poincaré characteristic gives

$$\dim H^1(G_S, U) - \dim H^2(G_S, U) = \dim H^0(G_S, U) + \dim U^*(1) = 0 + 2 = 2,$$

and substituting this as well as equation (34) into the sequence (31) we deduce the result. $\qquad \square$

**Remark 4.4.2** $H^1_f(K, U)$ is zero if and only if $H^1_f(K, A_U)$ has a trivial maximal divisible subgroup, and similarly for $H^1_f(K, A^*_U(1))$.

**Proposition 4.4.3**
$H^1_f(K, A_U)$ *is finite if and only if* $H^1_f(K, A^*_U(1))$ *is finite. in this case*

$$\#H^1_f(K, A^*_U(1)) = \#H^1_f(K, A_U) < \infty.$$

*Furthermore if* $k \geq 2$ *then both* $H^1_f(K, A_k)$ *and* $H^1_f(K, A^*_k(1))$ *are finite (see Section 5 for the situation for* $W_k$ *and* $W^*_k(1)$).

*Proof:* $H^1_f(K, A_U)$ and $H^1_f(K, A^*_U(1))$ are co-finitely generated. This follows from the fact that $H^1(G_S, U)$ has finite dimension.

Suppose $H^1_f(K, A_U)$ is finite. Then by Remark 4.4.2 $H^1_f(K, U) = 0$. By Proposition 4.4.1 $H^1_f(K, U^*(1)) = 0$, so $H^1_f(K, A^*_U(1))$ is also finite. This argument is reversible, and gives the first part of the Proposition.

On the other hand, if $k \geq 2$ then $H^1_f(K, A_k)$ is finite by Corollary 4.3.4 and Corollary 4.3.4, so by the first part $H^1_f(K, A^*_k(1))$ is also finite.

If both $H^1_f(K, A_U)$ and $H^1_f(K, A^*_U(1))$ are finite then we may utilise the generalised Cassels-Tate pairing of [Flach] which exhibits the two groups as duals of one another. The result is immediate. $\qquad\square$

**Proposition 4.4.4**
For $k \geq 1$ we have

$$\left( \frac{L(\overline{\psi}^k, k)}{\Omega^k} \right)_p = \left( \frac{(2\pi)^{k-1}}{(k-1)!} \frac{L(\psi^k, 1)}{\Omega^k} \right)_p .$$

*Proof:* This follows by substituting $s = k$ into the functional equation (7), and noting that since $p$ is inert in $K/\mathbb{Q}$, $\epsilon(\overline{\psi}^k, k)$ is prime to $p$. $\qquad\square$

**Corollary 4.4.5** *For $k \geq 1$ we have*

$$\#H^1_f(K, A^*_k(1)) = \left( \frac{(2\pi)^{2k-1}}{(k-1)!^2} \frac{L(\psi^k, 1)}{\Omega^k} \frac{L(\overline{\psi}^k, 1)}{\overline{\Omega}^k} \right)_p .$$

*If $k \geq 2$ then both sides are finite.*

*Proof:* This follows from Propositions 4.4.3 and 4.4.4. $\qquad\square$

## 4.5 Descent to $\mathbb{Q}$

In this section we indicate how our results may be descended to $\mathbb{Q}$. In other words we show how to prove analogous formulae over $\mathbb{Q}$, under the assumption (which we make from this point on) that our elliptic curve $E$ is defined over $\mathbb{Q}$.

A special case of a very general result is

**Lemma 4.5.1** *Let $c$ denote a generator for $G(K/\mathbb{Q})$. We have*

$$H^1_f(K, res(A_U)) = H^1_f(K, res(A_U))^{c=1} \oplus H^1_f(K, res(A_U))^{c=-1}.$$

*Here $res(A_U)$ denotes the restriction of the $G_\mathbb{Q}$-representation $A_U$ to $G_K$.*

*Proof:* This follows because the Selmer groups are $p$-torsion, and $p$ is odd. $\qquad\square$

**Lemma 4.5.2** *Suppose $H^1_f(K, res(A_U))$ has finite order. Then*

$$\#H^1_f(K, res(A_U))^{c=1} = \#H^1_f(K, res(A_U))^{c=-1}.$$

*Proof:* First note that $p \nmid d_k$. This is because $p$ divides $d_K$ if and only if $E$ has bad reduction at $p$, and we have assumed that $E$ has good reduction at $p$. The complex multiplication by $\sqrt{d_K}$ acts on $A_U$, hence by functoriality, and since $p \nmid d_K$ it induces an involution of $H^1_f(K, res(A_U))$ which swaps the two eigenspaces for $c$ (since $G_{\mathbb{Q}}$ acts $K$-semi-linearly). $\square$

**Lemma 4.5.3**
$$H^1_f(K, res(A_U))^{c=1} = H^1_f(\mathbb{Q}, A_U).$$

*Proof:* This follows easily from the fact that $p$ is odd, and $G(K/\mathbb{Q})$ has order 2. $\square$

**Remark 4.5.4** Notice that we also have

$$
\begin{aligned}
H^1_f(\mathbb{Q}, tw(A_U)) &= H^1_f(K, tw(A_U))^{c=1} \\
&= (tw(H^1_f(K, A_U)))^{c=1} \\
&= H^1_f(K, A_U)^{c=-1},
\end{aligned}
\tag{35}
$$

where $tw$ denotes a twist by the Galois character associated to $K$. Thus by Lemmas 4.5.2, 4.5.3, and by (35) we obtain

$$H^1_f(\mathbb{Q}, A_U) \cong H^1_f(\mathbb{Q}, tw(A_U)).$$

Let $\Omega_+$ denote the real period of $E$. In other words, if $\omega$ is the canonical differential attached to our choice of Weierstrass equation for $E$, then $\Omega_+$ is the integral of $\omega$ over the connected component of the identity on the real points of $E$. Equivalently, if $L$ is the lattice of periods of $E/\mathbb{C}$, $\Omega_+$ is the minimal positive real lattice point (it is well known that $L$ is fixed by complex conjugation, so such points exist). There is the following relation between $\Omega$ and $\Omega_+$.

**Lemma 4.5.5** $\Omega_+ = z\Omega$, where $z \in \mathcal{O}_K$ is prime to $6\sqrt{d_K}$.

*Proof:* Let $\Omega z$ be the real period of $L$. Then $\Omega z = \overline{\Omega} \overline{z}$. Hence $\Omega = \overline{\Omega} w$ for some $w \in K$. Taking norms we find $N_{K/\mathbb{Q}}(w) = 1$, so $w$ is a unit in $\mathcal{O}_K$. We have $\overline{z} = wz$. We consider cases.

If $w = 1$ then $z$ is real. By minimality $z = 1$, and $\Omega = \Omega_+$. If $w = -1$ then $z$ is purely imaginary. Hence by minimality $z = \sqrt{d_K}$.

If $w = \pm i$, we see that $K = \mathbb{Q}(i)$ and check (using minimality) that $z \mid 2$. Similarly, if $w^3 = 1$, then $K = \mathbb{Q}(\sqrt{-3})$, and $z \mid 3$. $\square$

**Theorem 4.5.6**
If $k \geq 1$ then

$$
\#H^1_f(\mathbb{Q}, A_k) = \left( \frac{L(\overline{\psi}^k, k)}{\Omega_+^k} \right)_p
\tag{36}
$$

and

$$
\#H^1_f(\mathbb{Q}, A_k^*(1)) = \left( \frac{L(\psi^k, k)}{\Omega_+^k} \right)_p .
\tag{37}
$$

Furthermore, if $k \geq 2$ then all expressions are finite.

*Proof:* By Lemmas 4.5.1, 4.5.2, and 4.5.3 we have that

$$\#H^1_f(K, res(A_k)) = \left(\#H^1_f(\mathbb{Q}, A_k)\right)^2.$$

Hence (36) follows by substituting these and (9) into Theorem 4.3.3. Note that after Lemma 4.5.5 and by our hypothesis on $p$ we may replace $\Omega$ by $\Omega_+$. The rest follows similarly, using Corollary 4.4.5. $\qquad\square$

# 5   Application to Chow Groups

Though results on the special values of $L$-functions are interesting in their own right, it turns out they they have applications in the study of cycle class maps. Langer and Raskind [LR] are able to show the finiteness of $CH_0(E \times E)\{p\}$, for $E/\mathbb{Q}$ an elliptic curve with complex multiplication by $\mathcal{O}_K$ having good reduction at $p > 3$, using Corollary 4.4.5 for the supersingular reduction case. In this section we show how our results may be used together with fundamental results of Nekovář [N2] to provide some information about higher cycle class maps for the products $E^d$, $d \geq 1$.

## 5.1   Results modulo torsion

Let $X$ be a smooth projective scheme over a number field $F$. There is a cycle class map into continuous étale cohomology (see [Jan]) for any prime $p$

$$cl_X : CH^i(X) \otimes \mathbb{Q}_p \longrightarrow H^{2i}(X_{et}, \mathbb{Q}_p(i)).$$

Writing $V^i = V^i(X) = H^i((X \otimes \overline{F})_{et}, \mathbb{Q}_p)$ one has the convergent Hoschild-Serre spectral sequence

$$E_2^{a,b} = H^a(F, V^b(X)(i)) \Rightarrow H^{a+b}(X_{et}, \mathbb{Q}_p(i))$$

inducing a filtration on its abutment. In fact it degenerates at the $E_2$ term, so

$$gr^j H^{2i}(X_{et}, \mathbb{Q}_p(i)) = H^j(F, V^{2i-j}(i)).$$

It is expected that there is a motivic version of the Hoschild-Serre spectral sequence yielding a filtration on rational Chow groups (motivic cohomology groups) strictly compatible with the above filtration on $l$-adic cohomology, and that the global dimension of the category of mixed motives over $F$ is 1. This would imply that

$$\mathrm{Im}\, cl_X \cap F^2 H^{2i}(X_{et}, \mathbb{Q}_p(i)) = 0. \tag{38}$$

One may attempt to prove (38) independently of the above motivic considerations, and we may make a first step by reducing it to a certain local-global principle.

Suppose $\mathcal{X}$ is a smooth projective model of $X$ over $\mathrm{Spec}\,\mathcal{O}_F[\frac{1}{N}]$ for some (square-free) integer $N$, and write

$$j : X \hookrightarrow \mathcal{X}$$

for the inclusion. Write $S = \{v \mid pN\}$. There are versions of the above constructions for $\mathcal{X}$ and one may show that

$$\text{Im}(cl_X) \cap F^2 H^{2i}(X_{et}, \mathbb{Q}_p(i)) = j^*(H),$$

where

$$H = \left( \text{Im}(cl_{\mathcal{X}}) \cap F^2 H^{2i}(\mathcal{X}_{et}, \mathbb{Q}_p(i)) \right).$$

The following is a special case of a fundamental theorem proved in [N2]. Fix $i$ and put $V = V^{2i-2}$, to ease notation.

**Theorem 5.1.1**
*Suppose that $X$ has potentially good reduction at every prime ideal of $F$. Then*

$$\text{Im}(cl_X) \cap F^2 H^{2n}(X_{et}, \mathbb{Q}_p(i))$$

*is a subquotient vector space of*

$$\mathcal{K}' = \ker \left( \alpha : H^2(G_{F,S}, V) \longrightarrow \bigoplus_{v \in S} H^2(F_v, V) \right).$$

On the other hand global Poitou-Tate duality gives $K'$ as the dual of

$$\mathcal{K} = \ker \left( \beta : H^1(G_{F,S}, V^*(1)) \longrightarrow \bigoplus_{v \in S} H^1(F_v, V^*(1)) \right),$$

and it is conjectured that $K$ always vanishes. While this local-global principle seems to be hard to prove in general, we may sometimes prove the triviality of

$$H^1_f(F, V^*(1)) = \ker \left( H^1(G_{F,S}, V^*(1)) \longrightarrow \bigoplus_{v \in S} \frac{H^1(F_v, V^*(1))}{H^1_f(F_v, V^*(1))} \right),$$

an a priori larger group. The dimension of $H^1_f(F, V^*(1))$ is conjecturally controlled by the order of zero of the $L$-function $L(V, s)$ at $s = 0$ (at least if $H^0(F, V^*(1)) = 0$). Of course $L(V, s)$ may vanish at $s = 0$, but there are interesing situations where one may show that it does not vanish, and in his case one might try to prove the part of the Bloch-Kato conjecture which tells us that $H^1_f(F, V^*(1))$ is trivial, hence yielding (38).

From now on let $F = K$, the quadratic imaginary field considered previously. Put $X/K = (E/K)^d$, the self product of $d$ copies of $E/K$, where $E/K$ is the elliptic curve with complex multiplication considered in previous chapters, and suppose for the moment that $E$ is defined over $\mathbb{Q}$. The Künneth formula gives $V$ as a direct sum of groups of representations of the form $V^1(E)^{\otimes_{\mathbb{Q}_p} 2i-2-2k}(-k)$ for varying $k$. The next step is thus to understand the representations $V^1(E)^{\otimes_{\mathbb{Q}_p} j}$ in terms of the $V_k$ considered in earlier sections. We work over $K$, since the results of (4.5) will yield results over $\mathbb{Q}$.

As a representation of $G_K$, some linear algebra yields

$$V^1(E)^{\otimes_{\mathbb{Q}_p} 2} = ((\Lambda^2_{\mathbb{Q}_p} V^1(E)) \otimes K_p) \oplus Sym^2_{K_p} V^1(E),$$

which by the Weil pairing is

$$K_p(-1) \oplus V^1(E)^{\otimes_{K_p} 2} = K_p(-1) \oplus V_2(-2).$$

Here we have identified $V^1(E)(1)$ with the rational Tate module $V_1$. Hence, again as $G_K$ representations, we find that $V^1(E)^{\otimes_{\mathbb{Q}_p} 2i-2-2k}(i-k)$ is a direct sum of representations of the form

$$W_m \cong V_2(-2)^{\otimes_{\kappa_p} m}(1+m),$$

for $0 \leq m \leq i-k-1$. Taking the Tate twisted dual, we find (using the Weil pairing here) that we need to bound the Selmer groups for the representations

$$W_m^*(1) \cong V_2(E)^{\otimes_{\kappa_p} m}(-m).$$

In general all these representations will occur with various (nonzero) multiplicities, but it is irrelevant to us what these multiplicities are.

In [Han], Han considers the Tate twisted representations

$$V_{j,k} = V_1(E)^{\otimes_{\kappa_p}(j+k)}(-j) \qquad j \geq 0 \qquad k > j+1,$$

whose Selmer group should conjecturally be related to the $L$-value

$$L(\overline{\psi}^{j+k}, k).$$

Here $E$ is as above. In the case $j = 0$ we find ourselves in the situation discussed in the first four sections of this article. Han is able to prove the following.

**Theorem 5.1.2**
Let $T_{j,k}$ be the canonical lattice in $V_{j,k}$. Let $p$ be a prime such that

- $p > 3$,

- $\psi^k \overline{\psi}^{-j}$ is nontrivial when restricted to the decomposition group of $K(E[p])/K$ at each prime above $p$.

Under these hypotheses, the Selmer group $H^1_f(K, V/T)$ is finite with order equal to

$$\left( (k-1)! \frac{(2\pi)}{\sqrt{d_K}} \frac{L(\overline{\psi}^{j+k}, k)}{\Omega_\infty^{j+k}} \right)_p.$$

**Remark 5.1.3** Note that here $p$ may be a prime of ordinary reduction, and possibly bad reduction for $E$. Actually Han goes further to show the Bloch-Kato conjecture holds for the corresponding motive. Also note that we have verified in our situation that the second hypothesis above always holds for $p > k$, (with $j = 0$, of course).

**Remark 5.1.4** If $p$ is prime to $6\mathfrak{f}$ then $\psi$ is surjective. If $p$ is supersingular then the condition on $\psi$ in the statement of the theorem will be valid for $p$ such that $p^2 - 1 \nmid n+1+p(n-1)$ for $0 \leq n \leq i-1$ (see Lemma 3.4.2). In particular this will hold for $p > 2d+1$.

If $p$ is ordinary and prime to $6\mathfrak{f}$, then the condition becomes $n+1+p(n-1) \equiv 0 \bmod p-1$, i.e $p-1 \nmid 2n$. Again this will hold for $p > 2d+1$.

Recall that we are interested in a representation of weight 0. Han's are of weight less than or equal to −2. On the other hand $V^*_{m-1,m+1}(1) = V_1(E)^{2m}(-m)$. We may thus use Proposition 4.4.1, and Han's result for $V_{m-1,m+1}$, together with the observation we have made before that $H^1_f(K, A^*_U(1))$ is finite if and only if $H^1_f(K, V^*_U(1))$ is trivial. We deduce the following, which is our main theorem modulo torsion. See the Appendix for remarks on the proof in the case of an ordinary prime.

**Theorem 5.1.5**
*Suppose $E$ is an elliptic curve over $\mathbb{Q}$, with good reduction at a prime $p > 3$. Let $i$ be a strictly positive integer. Assume that*

- *$p^2 - 1 \nmid n + 1 + p(n-1)$ if $E$ has supersingular reduction at $p$,*

- *$p - 1 \nmid 2n$ if $E$ has ordinary reduction at $p$,*

*for all $0 \le n \le i - 1$. Then*

$$\operatorname{Im}(cl_X(CH^i(X)) \cap F^2H^{2i}(X_{et}, \mathbb{Q}_p(i)) = 0,$$

*where $X$ is the d-fold product of $E$. The same is true for $E \otimes_{\mathbb{Q}} K$.*

$\square$

## 5.2 Integral results

Let our notation be as at the beginning of Section 5.1. Most of Section 5.1 carries through to the integral situation, given that the theorem of Han actually yields the finiteness of an integral Selmer group, not just the vanishing of $H^1_f(K, V)$ which is what was used in the last section. There are some points which require modification, however, and in fact our final theorem is stronger than the analgue of Theorem 5.1.5.

There is an integral cycle class map into continuous étale cohomology

$$cl_X : CH^i(X) \longrightarrow H^{2i}(X_{et}, \mathbb{Z}_p(i)).$$

Writing $T^i(X) = H^i((X \otimes \overline{F})_{et}, \mathbb{Z}_p)$ one has the Hoschild-Serre convergent spectral sequence

$$E_2^{a,b} = H^a(G_F, T^b(X)(i)) \Rightarrow H^{a+b}(X_{et}, \mathbb{Z}_p(i)).$$

This degenerates at the $E_2$ term up to finite groups, and in any case induces a filtration on its abutment. It is expected that

$$\operatorname{Im}(cl_X) \cap F^2H^{2i}(X_{et}, \mathbb{Z}_p(i))$$

is always finite. Suppose that $T^b$ is torsion free for all $b$. Then by the Weil conjectures $E_2^{0,b} = 0$ for all $b \ne 2i$. Since $cd_p(G_F) = 2$, we deduce

$$F^2(H^{2i}(X_{et}, \mathbb{Z}_p(i))) = gr^2(H^{2i}(X_{et}, \mathbb{Z}_p(i))) = H^2(G_F, T^{2i-2}(i)).$$

One may as before show that

$$\operatorname{Im}(cl_X) \cap F^2H^{2i}(X_{et}, \mathbb{Z}_p(i)) = j^*(H),$$

where
$$H = \left( \mathrm{Im}(cl_{\mathcal{X}}) \cap F^2 H^{2i}(\mathcal{X}_{et}, \mathbb{Z}_p(i)) \right).$$

On examining the local situation one finds that if $p > max\{2, 2i+1, \dim \mathcal{X}\}$, and if at every prime $v \in S$ $X_v$ achieves good reduction after base change to a field of degree prime to $p$ over $F_v$ then

$$\mathrm{Im}(cl_{X_v}) \cap F^2 H^{2i}((X_v)_{et}, \mathbb{Z}_p(i)) = 0$$

for every $v \in S$. There are two cases. If $v \mid p$ then this is Proposition III.5.5 in [N2]. If $v \nmid p$ then this is well known (see Proposition II.2.8 in [N2]).

**Remark 5.2.1** In fact it follows from III.5.6 in [N2] that we may also take $p = 2i+1$ for the choice of $X$ we shall shortly restrict to. However the conditions in Han's result require us to restrict to $p > 2i+1$.

Suppose from now on that we are in the case where the Hoschild - Serre spectral sequence degenerates (not just modulo finite groups). This will be the case later since all our cohomology groups will be torsion free. It follows from the compatibility of all the above constructions with base change that

$$H \hookrightarrow \ker \left( \alpha : H^2(G_{F,S}, T) \longrightarrow \bigoplus_{v \in S} H^2(F_v, T) \right),$$

where $T = T^{2i-2}(i)$. The Poitou-Tate duality sequence gives $\ker(\alpha)$ as the dual of

$$\mathcal{K} = \ker \left( \beta : H^1(G_{F,S}, T^*(1) \otimes \mathbb{Q}_p/\mathbb{Z}_p) \longrightarrow \bigoplus_{v \in S} H^1(F_v, T^*(1) \otimes \mathbb{Q}_p/\mathbb{Z}_p) \right).$$

Ideally one would like to prove that $\mathcal{K}$ was always finite, but this seems to be hard. Our approach is to estimate $\mathcal{K}$ by a Selmer group which we may be able to bound in terms of an $L$-value. Of course this does not work if the $L$-value is zero (but see [Ot]).

As before we put $F = K$, and $X/K = (E/K)^d$. The linear algebra considered in the last part works integrally, so long as we take $p > d$ (there are denominators involved in the various splittings of the representations).

The duality in Section 4.4 was dealt with on the integral level, so we may use it here, and similarly for the descent from $K$ to $\mathbb{Q}$. We have proved all but the last part of the following.

**Theorem 5.2.2**
Let $d$ be a positive integer, and let $p > 3$ be a rational prime such that

- $p > 2i+1$ (coming from the local result of Nekovář and the local condition of Han),

- $p > d$,

- $E$ has good reduction at $p$.

*Under these conditions*

$$\mathrm{Im}(cl_X(CH^i(X)) \cap F^2 H^{2i}(X_{et}, \mathbb{Z}_p(i))$$

*is finite, and is trivial for almost all primes. In fact*

$$\mathrm{Im}(cl_X(CH^i(X)))_{tors}$$

*is finite for $p$ as above, and trivial if $p$ also satisfies*

$$p \nmid (2\pi i)^j \frac{L(\psi^{2j}, -j)}{\Omega^{2j}}.$$

*Proof:* We need to prove only the final statement. We show the equality

$$H^{2i}(X_{et}, \mathbb{Z}_p(i))_{tors} = (F^2 H^{2i}(X_{et}, \mathbb{Z}_p(i)))_{tors}.$$

This together with the first statement will yield the final statement of the theorem. It is sufficient to show that both $gr^0 H^{2i}(X_{et}, \mathbb{Z}_p(i))$ and $gr^1 H^{2i}(X_{et}, \mathbb{Z}_p(i))$ are torsion free. But

$$gr^0 H^{2i}(X_{et}, \mathbb{Z}_p(i)) \cong (H^{2i}(\overline{X}_{et}, \mathbb{Z}_p(i)))^{G_\mathbb{Q}},$$

which is trivial by Künneth and because the Tate module of $E$ is torsion free. Furthermore,

$$gr^1 H^{2i}(X_{et}, \mathbb{Z}_p(i)) \cong H^1(\mathbb{Q}, H^{2i-1}(\overline{X}_{et}, \mathbb{Z}_p(i))),$$

and considering the long exact sequence corresponding to

$$0 \longrightarrow T \longrightarrow V \longrightarrow V/T \longrightarrow 0$$

we see that

$$H^1(\mathbb{Q}, H^{2i-1}(\overline{X}_{et}, \mathbb{Z}_p(i)))_{tors} \cong H^0(\mathbb{Q}, V/T),$$

where $T \cong H^{2i-1}(\overline{X}, \mathbb{Z}_p(i))$, and $V = T \otimes_{\mathbb{Z}_p} \mathbb{Q}_p$. the group $H^0(\mathbb{Q}, V/T)$ is trivial by the surjectivity of $\psi$, and our conditions on $p$.

## 5.3  Remarks on Ordinary Primes

We restricted to the case of supersingular reduction in the first four sections of this paper only because the case of ordinary reduction was dealt with by Li Guo in [LiGuo]. The proofs of Theorems 5.2.2 and 5.1.5 are, if anything, easier in the ordinary case. We give a few remarks on the differences in this situation for the interested reader.

If $E$ has good ordinary reduction at $p$ then $p$ splits in $K/\mathbb{Q}$,

$$(p) = \mathfrak{p}\overline{\mathfrak{p}}.$$

Write $K_\mathfrak{p}$ and $K_{\overline{\mathfrak{p}}}$ for the completion of $K$ at these primes. Both fields are isomorphic to $\mathbb{Q}_p$. Let $\pi$ be a uniformiser of $K_\mathfrak{p}$. The Tate module of $E$, is still free of rank 1 over $K \otimes \mathbb{Q}_p$, and this induces a decomposition of $G_K$ representations $V \cong V_\pi \oplus V_{\overline{\pi}}$. Being subrepresentations of $V$, both $V_\pi$ and $V_{\overline{\pi}}$ are crystalline when restricted to $G_{\mathbb{Q}_p}$.

In the Poitou-Tate sequence there are now two local conditions for $v \mid p$. We calculate that for the field $K_{\mathfrak{p}}$,

$$\dim D_{dR}^0(V_{\pi}) = 0, \qquad \dim D_{dR}^0(V_{\overline{\pi}}) = 1,$$

using the fact that $V_{\overline{\pi}}$ is unramified. By Proposition 3.1.2 this yields

$$\dim H_f^1(K_{\mathfrak{p}}, V_{\pi}) = 1, \qquad \dim H_f^1(K_{\overline{\mathfrak{p}}}, V_{\pi}) = 0.$$

The analogue of Proposition 4.4.1 follows immediately.

Everything else works as before: Han's result holds for ordinary primes. As already remarked, the bounding of the Selmer group in this case is originally due to Li Guo. The statement of Li Guo's result is the following (see [LiGuo]).

**Theorem 5.3.1**
Let $E$ be as before, and $p$ be a prime of good ordinary reduction for $E$. Let $V_{j,k}$ be as before, and suppose that $p > k + 1$, and $0 < j < k$. Then

$$\#(S^{str}(K, V_{j,k}) = \left( (2\pi)^j \frac{L(\psi^{j+k}, k)}{\Omega_{\infty}^{j+k}} \frac{L(\overline{\psi}^{j+k}, k)}{\overline{\Omega}_{\infty}^{j+k}} \right)_p.$$

Here $S^{str}(K, A_{j,k})$ denotes the strict Selmer group of Greenberg. Han uses the same Selmer groups as us (and Bloch-Kato), but see ([Flach]) for a comparison of these approaches.

# References

[BK] Bloch, S., Kato, K.: *L-functions and Tamagawa numbers of motives*, The Grothendieck Festschrift, Vol. 1, Birkhauser (1990), 333–400.

[CF] Cassels, J., Frohlich, A.:*Algebraic number theory*, Academic Press (1967).

[CS] Coates, J., Sydenham, A.: *On the symmetric square of a modular elliptic curve*, in 'Elliptic curves, modular forms, and Fermat's Last Theorem', Internat. Press (1995) 2–21.

[CW] Coates, J., Wiles, A.: *On the conjecture of Birch and Swinnerton-Dyer*, Invent. Math **39** (1977) 223-251.

[Flach] Flach, M.: *A generalisation of the Cassels-Tate pairing*, J. Reine. Angew. Math. **412** (1990) 113-127.

[Fon1] Fontaine, J-M.: *Le corps des périodes p-adiques*, Astérisque **223** (1994) 59–111.

[Fon2] ―――― *Représentations p-adiques semi-stables*, Astérisque **223** (1994) 113–184.

[Fon3] ―――― *Groupes p-divisible sur les corps locaux*, Astérisque 47–48 (1997)

[FPR] Fontaine, J-M., Perrin-Riou, B.: *Autour de conjecture de Bloch et Kato: cohomologie galoisienne et valeurs de fonctions L*, Motives (Seattle, WA, 1991), 599–706, Proc. Symp. Pure. Math., **55**, Part 1.

[Green] Greenberg, R.: *On the structure of certain Galois groups*, Invent. Math. **47** (1978) 85–99.

[Han] Han, B.: *On Bloch-Kato's Tamagawa number conjecture for Hecke characters of quadratic imaginary fields*, preprint available on the World Wide Web as *www.dpmms.cam.ac.uk/Algebraic − Number − Theory/0070/th1.dvi*.

[Hida] Hida, H.: *Elementary theory of L-functions and Eisenstein series*, Cambridge University Press (1993).

[Jan] Jannsen, U.: *Continuous étale cohomology*, Math. Ann. **280** (1988) 207–245.

[Kato] Kato, K.: *Lectures on the approach to Iwasawa Theory for Hasse-Weil L-functions via $B_{dR}$*, LNM 1553 (1993) 50–163.

[KM] Katz, N., Messing, W.: *Some consequences of the Riemann Hypothesis for varieties over finite fields*, Invent. Math. **23** (1974) 73–77.

[Lang] Lang, S.: *Algebraic Number Theory*, Springer-Verlag (1986).

[LR] Langer, A., Raskind, W.: *0-cycles on the self-product of a CM elliptic curve over $\mathbb{Q}$*, preprint.

[LiGuo] Guo, L.: *General Selmer groups and critical values of Hecke L–functions*, Math. Ann. **297** (1993), no.2, 221–233.

[N1] Nekovář, J.: *On p-adic height pairings, Séminaire des Théorie des Nombres*, Birkhauser (1993) 127–202.

[N2] —— *Syntomic cohomology and p-adic regulators*, preprint.

[Ot] Otsubo, N.: *Selmer groups and zero-cycles on the Fermat quartic surface*, preprint.

[PR] Perrin-Riou, B.: *Arithmétique des courbes elliptiques et théorie de Iwasawa*, Mém. Soc. Math. France **17** (1984)

[Rubin] Rubin, K.: *On the 'main conjecture' of Iwasawa theory for imaginary quadratic fields*, Invent. math., **103** (1991), 25–68.

[Rubin2] —— *Elliptic curves with complex multiplication and the conjecture of Birch and Swinnerton-Dyer*, preprint available on the World Wide Web as *math.stanford.edu/ rubin/cime.dvi*

[Schap] Schappacher, N.: *Periods of Hecke characters*, LNM **1301**, Springer Verlag, 1988.

[Serre] Serre, J-P.: *Corps Locaux*, Hermann, Paris, 1962.

[deShalit]  de Shalit, E.:*Iwasawa theory of elliptic curves with complex multiplication*, Academic Press (1987).

[Sil2]  Silverman, J.:*Advanced topics in the arithmetic of elliptic curves*, Springer-Verlag (1994).

[Tate]  Tate, J.: *p-divisible groups*, in *Proceedings of a Conference on Local Fields*, (1967), Springer, 158–183.

[Wiles]  Wiles, A.: *Higher explicit reciprocity laws*, Ann. Math. **107** (1978) 235–254.

[Win]  Wintenberger, J-P.:*Structure Galoisienne de limites projectives d'unités locales*, Compositio Math. **42** (1980/81) 89–103.